\documentclass[11pt,a4paper]{article}
\usepackage[utf8]{inputenc}
\usepackage[T1]{fontenc}
\usepackage[english]{babel}
\usepackage{lmodern}
\usepackage{amsmath,amssymb,amsfonts,amsthm,mathrsfs}
\usepackage{geometry}
\usepackage{hyperref}
\hypersetup{colorlinks=true,linkcolor=black,citecolor=blue,urlcolor=blue}
\newtheorem{theorem}{Theorem}[section]
\newtheorem{lemma}[theorem]{Lemma}
\newtheorem{proposition}[theorem]{Proposition}
\newtheorem{corollary}[theorem]{Corollary}
\newtheorem{example}[theorem]{Example}
\newtheorem{remark}[theorem]{Remark}
\newcommand{\tg}{\widetilde{g}}
\newcommand{\tS}{\widetilde{S}}
\newcommand{\tnabla}{\widetilde{\nabla}}
\newcommand{\tR}{\widetilde{R}}
\newcommand{\Lie}{\mathcal{L}}
\numberwithin{equation}{section}
\title{Rigidity of Ricci solitons for $\varphi$-related metallic semi-Riemannian structures on para-Kenmotsu manifolds}
\author{Yasemin KARABABA\\
\small Independent Researcher, Eski\c{s}ehir, T\"urkiye\\
\small \href{https://orcid.org/0000-0003-0485-9384}{ORCID: 0000-0003-0485-9384}}
\date{}
\begin{document}
\maketitle
\begin{abstract}
In this paper, we investigate Ricci solitons in the para-Kenmotsu setting with respect to a $\varphi$-related metallic semi-Riemannian metric arising from the underlying almost paracontact metric structure, under the assumption that the potential vector field is pointwise collinear with the Reeb vector field. We first establish the relation $\tS=S+2n\beta\bigl(g-\eta\otimes\eta\bigr)$ between the Ricci tensors of the original metric $g$ and the $\varphi$-related metallic metric $\tg$. We then show that the Ricci soliton equation implies that the original metric $g$ is of $\eta$-Einstein type. By combining this condition with the contracted Bianchi identity, we obtain a rigidity result implying $f=0$, $\lambda=\frac{2n}{C}$, $S=-2ng$, and consequently $\tS=-\frac{2n}{C}\tg$. Thus, every Ricci soliton in this setting whose potential vector field is pointwise collinear with the Reeb vector field is necessarily trivial, and both $g$ and $\tg$ are Einstein. Finally, we provide two five-dimensional para-Kenmotsu examples: one realizing the resulting trivial Einstein Ricci soliton and another showing that failure of the required Einstein condition prevents the existence of such a soliton.
\end{abstract}

\noindent\textbf{Keywords:} Para-Kenmotsu manifold, metallic structure,
$\varphi$-related semi-Riemannian metric, Ricci soliton,
Reeb vector field, $\eta$-Einstein manifold

\noindent\textbf{2020 Mathematics Subject Classification:} 53C24, 53C25, 53C50
	
	\section{Introduction}
	
	Almost paracontact metric geometry provides a natural framework for
	studying geometric structures on odd-dimensional semi-Riemannian
	manifolds. Among its distinguished subclasses, para-Kenmotsu manifolds
	are characterized by specific differential conditions involving the
	structure tensors and the Levi-Civita connection. Their curvature
	properties impose strong restrictions on the Reeb vector field and the
	Ricci tensor. The basic theory and curvature identities used in this
	paper follow Zamkovoy \cite{Zamkovoy2009,Zamkovoy2018}.
	
	Metallic structures have been studied from various geometric
	perspectives. Gezer and Karaman \cite{GezerKaraman2015} investigated
	metallic Riemannian structures, while \"Ozkan and Y{\i}lmaz
	\cite{OzkanYilmaz2018} studied metallic structures on differentiable
	manifolds. Metallic pseudo-Riemannian structures and their curvature
	properties were further investigated by Blaga and Nannicini
	\cite{BlagaNannicini2019,BlagaNannicini2020}. A $\varphi$-related metallic semi-Riemannian structure was constructed in
\cite{SolgunKarababa2025} from an almost paracontact metric structure, and the corresponding Levi-Civita connection and
	curvature tensor were investigated in the para-Kenmotsu and
	para-Sasakian cases. This construction provides the geometric setting considered in the present paper.
	
	Within this setting, it is natural to consider Ricci soliton equations with respect to the induced metallic metric. Ricci solitons and related soliton equations have previously been studied in
	para-Kenmotsu geometry. Blaga \cite{Blaga2015} investigated
	$\eta$-Ricci solitons and obtained quasi-Einstein and, under additional
	curvature assumptions, Einstein characterizations. Patra
	\cite{Patra2019} studied Ricci solitons and Ricci almost solitons on
	para-Kenmotsu manifolds and, in the Ricci almost soliton setting, also
	considered potential vector fields pointwise collinear with the Reeb
	vector field. More recently, various extensions such as $\eta$-Einstein solitons, conformal
$\eta$-Ricci solitons, and their rigidity characterizations on
para-Kenmotsu manifolds have attracted further interest; see, for example,
Yolda\c{s} \cite{Yoldas2023} and Dey and Uddin \cite{DeyUddin2022}.
In these studies, the soliton equation is imposed on a para-Kenmotsu metric. Here, we work with the induced metallic metric $\widetilde{g}$ and explicitly track its relation to the original metric $g$.
	
	The main purpose of this paper is to determine these restrictions when the potential
	vector field is pointwise collinear with the Reeb vector field. We
	first derive an explicit relation between the Ricci tensors of the
	original metric $g$ and the induced metallic metric $\widetilde{g}$.
	For a potential vector field of the form $V=f\xi$, we show that the
	soliton equation forces the Ricci tensor of the original metric $g$ to
	be of $\eta$-Einstein type. By combining the resulting identities with
	the contracted Bianchi identity, we obtain a rigidity theorem showing
	that the potential vector field necessarily vanishes and that both $g$
	and $\widetilde{g}$ are Einstein. Finally, we present two five-dimensional para-Kenmotsu examples. The first realizes the
	resulting trivial Einstein Ricci soliton, whereas the second shows that
	failure of the necessary Einstein condition obstructs the existence of
	a Reeb-directed Ricci soliton.
\section{Preliminaries}

\subsection{Almost paracontact metric and para-Kenmotsu structures}

Let $M^{2n+1}$ be a smooth manifold. The triple $(\varphi,\xi,\eta)$
defines an almost paracontact structure on $M$ if
\begin{equation}\label{eq:almost-paracontact}
	\varphi^2=I-\eta\otimes\xi,\qquad
	\eta(\xi)=1,\qquad
	\varphi\xi=0,\qquad
	\eta\circ\varphi=0.
\end{equation}
Assume further that $M$ is endowed with a semi-Riemannian metric $g$
satisfying
\begin{equation}\label{eq:paracontact-metric}
	g(\varphi X,\varphi Y)
	=
	-g(X,Y)+\eta(X)\eta(Y).
\end{equation}
Then $(M,\varphi,\xi,\eta,g)$ is called an almost paracontact metric
manifold \cite{Zamkovoy2009}. In particular,
\begin{equation}
	g(X,\xi)=\eta(X),\qquad
	g(\xi,\xi)=1,\qquad
	g(\varphi X,Y)=-g(X,\varphi Y).
\end{equation}

An almost paracontact metric manifold
$(M,\varphi,\xi,\eta,g)$ is said to be para-Kenmotsu provided that
\begin{equation}\label{eq:para-kenmotsu}
	(\nabla_X\varphi)Y
	=
	g(X,\varphi Y)\xi+\eta(Y)\varphi X,
\end{equation}
where $\nabla$ denotes the Levi-Civita connection of $g$
\cite{Zamkovoy2018}. It then follows that
\begin{align}
	\nabla_X\xi
	&=
	-X+\eta(X)\xi,
	\label{eq:nabla-xi}\\
	(\nabla_X\eta)(Y)
	&=
	-g(X,Y)+\eta(X)\eta(Y),
	\label{eq:nabla-eta}\\
	d\eta
	&=
	0.
	\label{eq:d-eta}
\end{align}

We adopt the curvature convention
\begin{equation}
	R(X,Y)Z
	=
	\nabla_X\nabla_YZ
	-\nabla_Y\nabla_XZ
	-\nabla_{[X,Y]}Z.
\end{equation}
With this convention, a para-Kenmotsu manifold satisfies
\begin{align}
	R(X,Y)\xi
	&=
	\eta(X)Y-\eta(Y)X,
	\\
	R(X,\xi)Y
	&=
	g(X,Y)\xi-\eta(Y)X,
	\\
	S(X,\xi)
	&=
	-2n\eta(X),
	\label{eq:para-kenmotsu-ricci}
\end{align}
where $S$ denotes the Ricci tensor of $g$ \cite{Zamkovoy2018}.
Moreover,
\begin{equation}\label{eq:lie-xi-g}
	\Lie_\xi g
	=
	-2g+2\eta\otimes\eta,
\end{equation}
where $\Lie_\xi$ denotes the Lie derivative with respect to the Reeb
vector field $\xi$.

\subsection{Metallic semi-Riemannian structures}

We recall some basic notions concerning metallic semi-Riemannian
structures; see \cite{BlagaNannicini2019,BlagaNannicini2020}.
Let $M$ be a smooth manifold. A $(1,1)$-tensor field $\Psi$ on $M$
is called a metallic structure if
\begin{equation}
	\Psi^2=p\Psi+qI,
\end{equation}
where $p$ and $q$ are real parameters. If, in addition, $M$ is endowed
with a semi-Riemannian metric $g$ such that $\Psi$ is $g$-symmetric,
namely,
\begin{equation}
	g(\Psi X,Y)=g(X,\Psi Y),
\end{equation}
then $(\Psi,g)$ is called a metallic semi-Riemannian structure and
$(M,\Psi,g)$ is called a metallic semi-Riemannian manifold
\cite{BlagaNannicini2019,BlagaNannicini2020}.

In the present paper, we restrict our attention to the parameter region
\begin{equation}
	p,q\in\mathbb{R},\qquad
	\Delta:=\sqrt{p^2+4q}>0,\qquad
	q\neq0,
\end{equation}
so that the characteristic roots are real and distinct. Set
\begin{equation*}\label{eq:metallic-roots}
	\sigma=\frac{p+\Delta}{2},
	\qquad
	\sigma^*=\frac{p-\Delta}{2}.
\end{equation*}
Then
\begin{equation}\label{eq:root-relations}
	\sigma+\sigma^*=p,\qquad
	\sigma\sigma^*=-q,\qquad
	\sigma-\sigma^*=\Delta.
\end{equation}

\section{The $\varphi$-related metallic semi-Riemannian structure}

A metallic structure associated with an almost paracontact structure was introduced in \cite{SolgunKarababa2025} by defining
\begin{equation}
	\Psi X
	=
	\sigma X+(p-2\sigma)\eta(X)\xi.
\end{equation}
Following the terminology introduced there, $\Psi$ is called the
$\varphi$-related metallic structure, or simply the related structure.
The corresponding semi-Riemannian metric is given by
\begin{equation}
	\tg(X,Y)
	=
	g(\Psi X,\Psi Y)+qg(X,Y).
\end{equation}
Equivalently,
\begin{equation}\label{eq:tilde-g}
	\tg
	=
	Ag+B\,\eta\otimes\eta,
\end{equation}
where
\begin{equation}\label{eq:ABC}
	A=p\sigma+2q,\qquad
	B=p(p-2\sigma),\qquad
	C:=A+B.
\end{equation}
Using \eqref{eq:root-relations}, these coefficients may be written as
\begin{equation*}\label{eq:ABC-roots}
	A=\sigma\Delta,\qquad
	B=-p\Delta,\qquad
	C=-\sigma^*\Delta.
\end{equation*}

Since $q\neq0$, both $\sigma$ and $\sigma^*$ are nonzero. Hence
$A\neq0$ and $C\neq0$, and therefore $\tg$ is nondegenerate.
Moreover,
\begin{equation}\label{eq:tilde-g-xi}
	\tg(X,\xi)=C\eta(X),
	\qquad
	\tg(\xi,\xi)=C.
\end{equation}

For later use, set
\begin{equation}
	\beta:=\frac{B}{C}.
\end{equation}
Then
\begin{equation}\label{eq:beta-relations}
	\beta=\frac{p}{\sigma^*},
	\qquad
	\beta=1-\frac{A}{C}.
\end{equation}

In the para-Kenmotsu case, the Levi-Civita connection and curvature
relations obtained in \cite{SolgunKarababa2025} are
\begin{equation*}\label{eq:connection-transform}
	\tnabla_XY
	=
	\nabla_XY+\beta g(\varphi X,\varphi Y)\xi,
\end{equation*}
and
\begin{equation}\label{eq:curvature-transform}
	\tR(X,Y)Z
	=
	R(X,Y)Z
	+\beta\bigl\{
	g(\varphi X,\varphi Z)Y
	-g(\varphi Y,\varphi Z)X
	\bigr\}.
\end{equation}
respectively. In particular,
\begin{align}
	\tR(X,Y)\xi
	&=
	R(X,Y)\xi,
	\\
	\tR(X,\xi)Y
	&=
	R(X,\xi)Y
	+\beta g(\varphi X,\varphi Y)\xi.
\end{align}

The preceding curvature relation allows us to express the Ricci tensor of the induced metallic metric $\tg$ in terms of the Ricci tensor of the original para-Kenmotsu metric $g$.

\begin{proposition}\label{prop:ricci-transform}
Let $(M^{2n+1},\varphi,\xi,\eta,g)$ be a para-Kenmotsu manifold, and let $\tg$ be the $\varphi$-related metallic semi-Riemannian metric defined by \eqref{eq:tilde-g}. If $S$ and $\tS$ denote the Ricci tensors of $g$ and $\tg$, respectively, then
\begin{equation}\label{eq:ricci-transform}
\tS(Y,Z)=S(Y,Z)+2n\beta\bigl\{g(Y,Z)-\eta(Y)\eta(Z)\bigr\}
\end{equation}
for all $Y,Z\in\Gamma(TM)$.
\end{proposition}

\begin{proof}
	For fixed $Y$ and $Z$, let
	$\{E_1,\ldots,E_{2n+1}\}$ be a local $g$-pseudo-orthonormal frame
	such that
	$g(E_i,E_i)=\varepsilon_i\in\{-1,1\}$.
	Since the Ricci tensor is the trace of the endomorphism $X\mapsto\tR(X,Y)Z$, this trace may be computed with respect to the $g$-pseudo-orthonormal frame. Hence
	\begin{equation}
		\tS(Y,Z)
		=
		\sum_{i=1}^{2n+1}
		\varepsilon_i
		g\bigl(\tR(E_i,Y)Z,E_i\bigr).
	\end{equation}
	Using \eqref{eq:curvature-transform}, we obtain
	\begin{align}
		\tS(Y,Z)
		=
		S(Y,Z)
		+\beta
		\sum_{i=1}^{2n+1}
		\varepsilon_i
		g(\varphi E_i,\varphi Z)g(Y,E_i)
		-\beta
		g(\varphi Y,\varphi Z)
		\sum_{i=1}^{2n+1}
		\varepsilon_i g(E_i,E_i).
	\end{align}
	By pseudo-orthonormal completeness,
	\begin{equation}
		\sum_{i=1}^{2n+1}
		\varepsilon_i
		g(\varphi E_i,\varphi Z)g(Y,E_i)
		=
		g(\varphi Y,\varphi Z),
	\end{equation}
	while
	\begin{equation}
		\sum_{i=1}^{2n+1}
		\varepsilon_i g(E_i,E_i)
		=
		2n+1.
	\end{equation}
	Therefore,
	\begin{equation}
		\tS(Y,Z)
		=
		S(Y,Z)-2n\beta g(\varphi Y,\varphi Z).
	\end{equation}
	Finally, using \eqref{eq:paracontact-metric}, we obtain
	\begin{equation}
		\tS(Y,Z)
		=
		S(Y,Z)
		+2n\beta
		\bigl\{
		g(Y,Z)-\eta(Y)\eta(Z)
		\bigr\},
	\end{equation}
	which proves \eqref{eq:ricci-transform}.
\end{proof}

As an immediate consequence, the Ricci tensor of the induced metric has the same component in the Reeb direction as that of the original para-Kenmotsu metric.

\begin{corollary}\label{cor:ricci-xi}
	For every $X\in\Gamma(TM)$,
	\begin{equation}\label{eq:tilde-ricci-xi}
		\tS(X,\xi)=-2n\eta(X),
		\qquad
		\tS(\xi,\xi)=-2n.
	\end{equation}
\end{corollary}

\begin{proof}
Setting $Z=\xi$ in \eqref{eq:ricci-transform} and using $g(X,\xi)=\eta(X)$ and $\eta(\xi)=1$, we obtain $\tS(X,\xi)=S(X,\xi)$. Hence, by \eqref{eq:para-kenmotsu-ricci}, $\tS(X,\xi)=-2n\eta(X)$. Taking $X=\xi$ and using $\eta(\xi)=1$ gives $\tS(\xi,\xi)=-2n$.
\end{proof}

\section{Ricci solitons with a Reeb-directed potential field}

We use the following sign convention for the Ricci soliton equation:
\begin{equation}\label{eq:soliton}
	\frac{1}{2}\Lie_V\tg+\tS+\lambda\tg=0,
\end{equation}
A triple $(\tg,V,\lambda)$ satisfying \eqref{eq:soliton} is called a Ricci soliton, where $V$ is the potential vector field and $\lambda\in\mathbb{R}$ is a constant \cite{Hamilton1988}.
Throughout this section, let $(M^{2n+1},\varphi,\xi,\eta,g)$ be a para-Kenmotsu manifold, and let $\tg$ denote the $\varphi$-related metallic semi-Riemannian metric defined in Section 3. The Ricci soliton equation is considered with respect to $\tg$.
We assume that the potential vector field is pointwise collinear with the Reeb vector field, namely,
\begin{equation*}\label{eq:potential}
V=f\xi,
\end{equation*}
where $f\in C^\infty(M)$.
We first determine the Lie derivative of the $\varphi$-related metallic metric along the Reeb vector field.

\begin{lemma}\label{lem:lie-xi}
Let $(M^{2n+1},\varphi,\xi,\eta,g)$ be a para-Kenmotsu manifold, and let $\tg$ be the $\varphi$-related metallic semi-Riemannian metric defined by \eqref{eq:tilde-g}. Then the Lie derivative of $\tg$ along the Reeb vector field $\xi$ satisfies
	\begin{equation*}\label{eq:lie-xi-tildeg}
\Lie_\xi\tg=-2\tg+2C\,\eta\otimes\eta,
\end{equation*}
where $C=A+B$ is given by \eqref{eq:ABC}.
\end{lemma}

\begin{proof}
	Since $A$ and $B$ are constant, it follows from \eqref{eq:tilde-g} that
	\begin{equation*}\label{eq:lie-linear}
\Lie_\xi\tg=A\Lie_\xi g+B\Lie_\xi(\eta\otimes\eta).
\end{equation*}
	For a para-Kenmotsu manifold,
	$\Lie_\xi(\eta\otimes\eta)=0$ \cite{Blaga2015}.
	Using \eqref{eq:lie-xi-g}, we obtain
	\begin{equation*}
		\begin{aligned}
			\Lie_\xi\tg
			&=-2Ag+2A\,\eta\otimes\eta\\
			&=-2\tg+2(A+B)\eta\otimes\eta\\
			&=-2\tg+2C\,\eta\otimes\eta.
		\end{aligned}
	\end{equation*}
\end{proof}

The preceding identity extends to potential vector fields of the form $V=f\xi$.

\begin{lemma}\label{lem:lie-fxi}
	Let $(M^{2n+1},\varphi,\xi,\eta,g)$ be a para-Kenmotsu manifold, let $\tg$ be the $\varphi$-related metallic semi-Riemannian metric, and let $V=f\xi$, where $f\in C^\infty(M)$. Then
	\begin{equation}\label{eq:lie-fxi}
		\Lie_{f\xi}\tg
		=
		-2f\tg
		+2fC\,\eta\otimes\eta
		+C\bigl(df\otimes\eta+\eta\otimes df\bigr).
	\end{equation}
\end{lemma}

\begin{proof}
	For any semi-Riemannian metric $h$ and vector field $X$, one has
	\begin{equation}
		\Lie_{fX}h
		=
		f\Lie_Xh
		+df\otimes X^{\flat_h}
		+X^{\flat_h}\otimes df,
	\end{equation}
	where $X^{\flat_h}=h(X,\cdot)$.
	By \eqref{eq:tilde-g-xi},
	$\xi^{\flat_{\tg}}=\tg(\xi,\cdot)=C\eta$.
	Therefore, \eqref{eq:lie-fxi} follows directly from
	Lemma \ref{lem:lie-xi}.
\end{proof}

We now substitute the preceding expression into the Ricci soliton equation.

\begin{proposition}\label{prop:tS-soliton}
	Let $(M^{2n+1},\varphi,\xi,\eta,g)$ be a para-Kenmotsu manifold equipped with the $\varphi$-related metallic semi-Riemannian metric $\tg$. If $(\tg,V,\lambda)$ defines a Ricci soliton with $V=f\xi$, then
	\begin{equation}\label{eq:tS-soliton}
		\tS
		=
		(f-\lambda)\tg
		-fC\,\eta\otimes\eta
		-\frac{C}{2}\bigl(df\otimes\eta+\eta\otimes df\bigr).
	\end{equation}
\end{proposition}

\begin{proof}
	Substituting \eqref{eq:lie-fxi} into the Ricci soliton equation
	\eqref{eq:soliton} and solving for $\tS$ yields
	\eqref{eq:tS-soliton}.
\end{proof}

The Reeb component of \eqref{eq:tS-soliton} determines the differential of the function $f$.

\begin{theorem}\label{thm:df}
	Let $(M^{2n+1},\varphi,\xi,\eta,g)$ be a para-Kenmotsu manifold equipped with the $\varphi$-related metallic semi-Riemannian metric $\tg$. If $(\tg,V,\lambda)$ is a Ricci soliton with $V=f\xi$, then
	\begin{equation}\label{eq:df}
		df
		=
		\left(\frac{2n}{C}-\lambda\right)\eta.
	\end{equation}
	In particular,
	\begin{equation}\label{eq:f-derivatives}
		X(f)=0
		\quad\text{for every }X\in\ker\eta,
		\qquad
		\xi(f)=\frac{2n}{C}-\lambda.
	\end{equation}
\end{theorem}

\begin{proof}
	Evaluating \eqref{eq:tS-soliton} at $(X,\xi)$ and using
	\eqref{eq:tilde-ricci-xi} and \eqref{eq:tilde-g-xi}, we obtain
	\begin{equation}
		-2n\eta(X)
		=
		-\lambda C\eta(X)
		-\frac{C}{2}df(X)
		-\frac{C}{2}\xi(f)\eta(X).
	\end{equation}
	Since $C\neq0$, it follows that
	\begin{equation}
		df(X)
		=
		\left(
		\frac{4n}{C}-2\lambda-\xi(f)
		\right)\eta(X).
	\end{equation}
	Setting $X=\xi$ gives
\begin{equation}\label{eq:xi-f-intermediate}
\xi(f)=\frac{4n}{C}-2\lambda-\xi(f),
\end{equation}
and hence
\begin{equation*}\label{eq:xi-f}
\xi(f)=\frac{2n}{C}-\lambda.
\end{equation*}
	Substituting this relation into the preceding equation yields
	\eqref{eq:df}. The identities in \eqref{eq:f-derivatives} follow immediately.
\end{proof}

We next derive the consequences of \eqref{eq:df} for the Ricci tensor of the original para-Kenmotsu metric $g$. By Theorem \ref{thm:df},
\begin{equation}
	df\otimes\eta+\eta\otimes df
	=
	2\left(\frac{2n}{C}-\lambda\right)\eta\otimes\eta.
\end{equation}
Substituting this relation into \eqref{eq:tS-soliton}, we obtain
\begin{equation}\label{eq:tS-eta}
	\tS
	=
	(f-\lambda)\tg
	+\bigl\{C(\lambda-f)-2n\bigr\}\eta\otimes\eta.
\end{equation}

Combining \eqref{eq:tS-eta} with the Ricci transformation formula \eqref{eq:ricci-transform} gives the following necessary condition on the Ricci tensor of the original metric.

\begin{proposition}\label{prop:eta-einstein}
Let $(M^{2n+1},\varphi,\xi,\eta,g)$ be a para-Kenmotsu manifold equipped with the $\varphi$-related metallic semi-Riemannian metric $\tg$. If $\tg$ admits a Ricci soliton $(\tg,V,\lambda)$ with $V=f\xi$, then the Ricci tensor $S$ of the original metric $g$ is necessarily of $\eta$-Einstein type:
	\begin{equation}\label{eq:S-explicit}
		S=ag+b\,\eta\otimes\eta,
	\end{equation}
	where
	\begin{align}
		a&=A(f-\lambda)-2n\beta,
		\label{eq:a}\\
		b&=-A(f-\lambda)-2n+2n\beta.
		\label{eq:b}
	\end{align}
\end{proposition}

\begin{proof}
	Substituting \eqref{eq:tilde-g} into \eqref{eq:tS-eta}, we obtain
	\begin{equation}
		\tS
		=
		A(f-\lambda)g
		+
		\bigl\{
		B(f-\lambda)+C(\lambda-f)-2n
		\bigr\}\eta\otimes\eta.
	\end{equation}
	On the other hand, \eqref{eq:ricci-transform} gives
	\begin{equation}
		\tS
		=
		S+2n\beta g-2n\beta\,\eta\otimes\eta.
	\end{equation}
	Comparing the coefficients of $g$ and $\eta\otimes\eta$ yields
	\eqref{eq:a} and \eqref{eq:b}. In particular,
	\begin{equation}\label{eq:a+b}
		a+b=-2n,
	\end{equation}
	which is consistent with \eqref{eq:para-kenmotsu-ricci}.
\end{proof}

The next result records the scalar-curvature relation that will be used in the proof of the main rigidity theorem.

\begin{proposition}\label{prop:scalar-bianchi}
	Let $(M^{2n+1},\varphi,\xi,\eta,g)$ be a para-Kenmotsu manifold equipped with the $\varphi$-related metallic semi-Riemannian metric $\tg$. Suppose that $(\tg,V,\lambda)$ is a Ricci soliton with $V=f\xi$.
	Then the scalar curvature $r$ of the original metric $g$ satisfies
	\begin{equation}\label{eq:r}
		r
		=
		2nA(f-\lambda)-2n-4n^2\beta,
	\end{equation}
	and
	\begin{equation}\label{eq:xi-r}
		\xi(r)
		=
		2\bigl\{r+2n(2n+1)\bigr\}.
	\end{equation}
\end{proposition}

\begin{proof}
	Taking the $g$-trace of \eqref{eq:S-explicit} gives
	$r=(2n+1)a+b$. Using \eqref{eq:a+b}, this reduces to
	$r=2na-2n$. Substitution of \eqref{eq:a} then yields
	\eqref{eq:r}.
	
	Moreover, \eqref{eq:a+b} together with
	$r=(2n+1)a+b$ gives
	\begin{equation*}
		a
		=
		1+\frac{r}{2n},
		\qquad
		b
		=
		-\left(2n+1+\frac{r}{2n}\right).
	\end{equation*}
	
	We now use the contracted Bianchi identity
	\begin{equation}\label{eq:contracted-bianchi}
		\operatorname{div}S=\frac{1}{2}dr.
	\end{equation}
	Since $S=ag+b\,\eta\otimes\eta$ and $\nabla g=0$, we have
	$\operatorname{div}(ag)=da$.
	
	Furthermore, $\operatorname{div}(b\,\eta\otimes\eta)=\xi(b)\eta+b(\operatorname{div}\xi)\eta+b\nabla_\xi\eta$. By \eqref{eq:nabla-xi} and \eqref{eq:nabla-eta}, $\operatorname{div}\xi=-2n$, $\nabla_\xi\eta=0$.
Therefore,
	\begin{equation*}
		\operatorname{div}(b\,\eta\otimes\eta)
		=
		\bigl\{\xi(b)-2nb\bigr\}\eta.
	\end{equation*}
	
	Since
	$da=\frac{1}{2n}dr$, $\xi(b)=-\frac{1}{2n}\xi(r)$, we obtain
	\begin{equation}\label{eq:divS}
		\operatorname{div}S
		=
		\frac{1}{2n}dr
		-\frac{1}{2n}\xi(r)\eta
		-2nb\,\eta.
	\end{equation}
	Evaluating \eqref{eq:contracted-bianchi} and \eqref{eq:divS} on $\xi$
	gives
\begin{equation*}\label{eq:bianchi-reeb}
\frac{1}{2}\xi(r)=-2nb.
\end{equation*}
	Substituting the above expression for $b$ yields
	\eqref{eq:xi-r}.
\end{proof}

\section{Rigidity of Reeb-directed Ricci solitons}
We are now in a position to establish the main rigidity result of the paper.

\begin{theorem}[Main rigidity theorem]\label{thm:main}
Let $(M^{2n+1},\varphi,\xi,\eta,g)$ be a para-Kenmotsu manifold, and let $(\Psi,\tg)$ be the $\varphi$-related metallic semi-Riemannian structure. Assume that $p,q\in\mathbb{R}$, $\Delta=\sqrt{p^2+4q}>0$, $q\neq0$. If $(\tg,V,\lambda)$ is a Ricci soliton whose potential vector field is pointwise collinear with the Reeb vector field, then the Ricci soliton is trivial. More precisely,
\begin{equation*}\label{eq:main-values}
f=0,\qquad V=0,\qquad\lambda=\frac{2n}{C},\qquad r=-2n(2n+1),
\end{equation*}
where $r$ denotes the scalar curvature of the original metric $g$. Moreover,
\begin{equation*}\label{eq:main-einstein}
S=-2ng,\qquad\tS=-\frac{2n}{C}\tg.
\end{equation*}
In particular, both $g$ and $\tg$ are Einstein metrics.
\end{theorem}

\begin{proof}
By Theorem \ref{thm:df}, $\xi(f)=\frac{2n}{C}-\lambda$. Since $C$ and $\lambda$ are constant, $\xi(f)$ is constant. Differentiating \eqref{eq:r} along the Reeb vector field $\xi$ gives
\begin{equation}\label{eq:xi-r-const}
\xi(r)=2nA\,\xi(f)=2nA\left(\frac{2n}{C}-\lambda\right).
\end{equation}
Thus $\xi(r)$ is constant, and consequently
\begin{equation}\label{eq:xi-xi-r-zero}
\xi\xi(r)=0.
\end{equation}
On the other hand, differentiating \eqref{eq:xi-r} along $\xi$ yields $\xi\xi(r)=2\xi(r)$. Combining this relation with \eqref{eq:xi-xi-r-zero}, we obtain $\xi(r)=0$. Substituting this into \eqref{eq:xi-r} gives
\begin{equation}\label{eq:scalar-rigid}
r=-2n(2n+1).
\end{equation}
Using \eqref{eq:xi-r-const}, we then have $2nA\left(\frac{2n}{C}-\lambda\right)=0$. Since $n\geq1$ and $A\neq0$, it follows that
\begin{equation}\label{eq:lambda-rigid}
\lambda=\frac{2n}{C}.
\end{equation}
Hence, by \eqref{eq:df}, $df=0$, so that $f$ is constant.
We next determine the Ricci tensor of the original metric. Taking the $g$-trace of \eqref{eq:S-explicit} gives $r=(2n+1)a+b$. Together with \eqref{eq:a+b}, this may be written as $r=2na-2n$. Using \eqref{eq:scalar-rigid}, we obtain $a=-2n$, $b=0$. Therefore,
\begin{equation*}\label{eq:S-rigid}
S=-2ng.
\end{equation*}
It remains to determine the constant value of $f$. From \eqref{eq:b}, together with $b=0$, \eqref{eq:lambda-rigid}, and \eqref{eq:beta-relations}, we obtain
\begin{equation*}\label{eq:f-zero-calculation}
\begin{aligned}
0&=-A\left(f-\frac{2n}{C}\right)-2n+2n\left(1-\frac{A}{C}\right)\\
 &=-Af.
\end{aligned}
\end{equation*}
Since $A\neq0$, it follows that $f=0$. Consequently, $V=f\xi=0$.
Finally, since $V=0$, the Ricci soliton equation \eqref{eq:soliton} reduces to
\begin{equation*}\label{eq:soliton-zero}
\tS+\lambda\tg=0.
\end{equation*}
Using \eqref{eq:lambda-rigid}, we obtain
\begin{equation*}\label{eq:tS-rigid}
\tS=-\frac{2n}{C}\tg.
\end{equation*}
Hence both $g$ and $\tg$ are Einstein.
\end{proof}

The preceding theorem immediately rules out nontrivial Ricci solitons whose potential vector fields are directed along the Reeb vector field.
\begin{corollary}\label{cor:no-nontrivial-soliton}
Let $(M^{2n+1},\varphi,\xi,\eta,g)$ be a para-Kenmotsu manifold equipped with the $\varphi$-related metallic semi-Riemannian metric $\tg$. Then $\tg$ admits no nontrivial Ricci soliton whose potential vector field is pointwise collinear with the Reeb vector field.
\end{corollary}

\begin{remark}
No additional sign condition on $C$ is required in Theorem \ref{thm:main}. Since
\begin{equation*}
C=-\sigma^*\Delta\neq0
\end{equation*}
and
\begin{equation*}
\lambda=\frac{2n}{C},
\end{equation*}
the sign of the soliton constant $\lambda$ is determined by the admissible metallic parameters through the sign of $C$.
\end{remark}

\section{Examples}
The following two five-dimensional examples illustrate two complementary aspects of Theorem \ref{thm:main}. Both examples are para-Kenmotsu manifolds. The first realizes the trivial Einstein Ricci soliton obtained in Theorem \ref{thm:main}, whereas the second shows that failure of the necessary Einstein condition prevents the existence of a Ricci soliton whose potential vector field is pointwise collinear with the Reeb vector field. Thus, the distinction between the two examples is geometric rather than dimensional.

\subsection{A five-dimensional realization of the trivial soliton}
\begin{example}
Let $M^5=\mathbb{R}^5(t,u,v,x,y)$ be endowed with the semi-Riemannian metric $g=dt^2+2e^{-2t}\,du\,dv+2e^{-2t}\,dx\,dy$. Set
\begin{equation*}
\xi=\partial_t,\qquad\eta=dt,
\end{equation*}
and define the $(1,1)$-tensor field $\varphi$ by
\begin{equation*}\label{eq:five-realization-phi}
\varphi\partial_u=\partial_u,\quad\varphi\partial_v=-\partial_v,\quad\varphi\partial_x=\partial_x,\quad\varphi\partial_y=-\partial_y,\quad\varphi\xi=0.
\end{equation*}
Then \eqref{eq:almost-paracontact} and \eqref{eq:paracontact-metric} are satisfied.
The nonzero covariant derivatives of the coordinate vector fields, up to the symmetry of the Levi-Civita connection, are
\begin{gather*}
\nabla_{\partial_u}\partial_v=e^{-2t}\xi,\quad\nabla_{\partial_x}\partial_y=e^{-2t}\xi,\quad\nabla_{\partial_u}\xi=-\partial_u,\\
\nabla_{\partial_v}\xi=-\partial_v,\quad\nabla_{\partial_x}\xi=-\partial_x,\quad\nabla_{\partial_y}\xi=-\partial_y.
\end{gather*}
In particular, \eqref{eq:nabla-xi} holds, and a direct verification on the coordinate frame gives \eqref{eq:para-kenmotsu}.
Hence $(M,\varphi,\xi,\eta,g)$ is a para-Kenmotsu manifold.
A direct curvature computation gives
\begin{equation*}\label{eq:five-realization-ricci}
S=-4g,\qquad r=-20.
\end{equation*}
Since $\dim M=5$, we have $n=2$. Theorem \ref{thm:main} therefore prescribes $\lambda=\frac{4}{C}$. We now verify directly that this value, together with $V=0$, satisfies the Ricci soliton equation.
By Proposition \ref{prop:ricci-transform},
\begin{align*}
\tS&=S+4\beta\bigl(g-\eta\otimes\eta\bigr)\\
&=-4g+4\beta\bigl(g-\eta\otimes\eta\bigr)\\
&=4(\beta-1)g-4\beta\,\eta\otimes\eta\\
&=-\frac{4A}{C}g-\frac{4B}{C}\eta\otimes\eta\\
&=-\frac{4}{C}\tg,
\end{align*}
where \eqref{eq:beta-relations} and \eqref{eq:tilde-g} have been used. Consequently, $\tS+\frac{4}{C}\tg=0$.
Since $V=0$, the Lie derivative term in \eqref{eq:soliton} vanishes. Thus
\begin{equation*}
V=0,\qquad\lambda=\frac{4}{C}
\end{equation*}
satisfy the Ricci soliton equation.
Moreover, by \eqref{eq:tilde-g} and \eqref{eq:tilde-g-xi},
\begin{equation}\label{eq:five-realization-restriction}
\tg|_{\ker\eta}=A\,g|_{\ker\eta},\qquad\tg(\xi,\xi)=C.
\end{equation}
Since $A\neq0$ and $C\neq0$, the induced metallic metric $\tg$ is nondegenerate.
Therefore, this five-dimensional para-Kenmotsu manifold realizes the trivial Einstein Ricci soliton described in Theorem \ref{thm:main}.
\end{example}

\subsection{A five-dimensional obstruction example}
\begin{example}
Let $M^5=\mathbb{R}^5(t,u,v,x,y)$ be endowed with the semi-Riemannian metric $g=dt^2+2e^{-2t}\,du\,dv+2e^{-2t+2xy}\,dx\,dy$.
Set
\begin{equation*}
\xi=\partial_t,\qquad\eta=dt,
\end{equation*}
and define $\varphi$ by
\begin{equation}\label{eq:five-obstruction-phi}
\varphi\partial_u=\partial_u,\quad\varphi\partial_v=-\partial_v,\quad\varphi\partial_x=\partial_x,\quad\varphi\partial_y=-\partial_y,\quad\varphi\xi=0.
\end{equation}
The relevant nonzero covariant derivatives of the coordinate vector fields are
\begin{gather*}
\nabla_{\partial_u}\partial_v=e^{-2t}\xi,\quad\nabla_{\partial_x}\partial_y=e^{-2t+2xy}\xi,\quad\nabla_{\partial_u}\xi=-\partial_u,\quad\nabla_{\partial_v}\xi=-\partial_v,\\
\nabla_{\partial_x}\xi=-\partial_x,\quad\nabla_{\partial_y}\xi=-\partial_y,\quad\nabla_{\partial_x}\partial_x=2y\partial_x,\quad\nabla_{\partial_y}\partial_y=2x\partial_y.
\end{gather*}
These relations show that \eqref{eq:nabla-xi} is satisfied. Moreover, a direct verification with respect to the coordinate frame confirms \eqref{eq:para-kenmotsu}. Hence, $(M,\varphi,\xi,\eta,g)$ is again a para-Kenmotsu manifold.
The nonzero Ricci components needed below are
\begin{equation*}\label{eq:five-obstruction-components}
S(\xi,\xi)=-4,\qquad S(\partial_u,\partial_v)=-4e^{-2t},
\end{equation*}
and
\begin{equation*}\label{eq:five-obstruction-ricci}
S(\partial_x,\partial_y)=-4e^{-2t+2xy}-2.
\end{equation*}
On the other hand,
\begin{equation*}\label{eq:five-obstruction-comparison}
-4g(\partial_x,\partial_y)=-4e^{-2t+2xy}.
\end{equation*}
Therefore,
\begin{equation}\label{eq:five-not-einstein}
S\neq-4g.
\end{equation}
Since $\dim M=5$, we have $n=2$. Suppose that the induced metallic metric $\tg$ admits a Ricci soliton whose potential vector field is of the form
\begin{equation*}
V=f\xi.
\end{equation*}
Then Theorem \ref{thm:main} would necessarily imply
\begin{equation}\label{eq:five-obstruction-necessary}
S=-4g,
\end{equation}
which contradicts \eqref{eq:five-not-einstein}. Hence the induced metallic metric $\tg$ admits no Ricci soliton whose potential vector field is pointwise collinear with the Reeb vector field.
Finally, by \eqref{eq:tilde-g} and \eqref{eq:tilde-g-xi},
\begin{equation*}\label{eq:five-obstruction-restriction}
\tg|_{\ker\eta}=A\,g|_{\ker\eta},\qquad\tg(\xi,\xi)=C.
\end{equation*}
Since $A\neq0$ and $C\neq0$, the induced metallic metric $\tg$ is nondegenerate.
Thus, this example satisfies the para-Kenmotsu condition but fails the necessary Einstein condition in Theorem \ref{thm:main}, and therefore cannot admit a Reeb-directed Ricci soliton.
\end{example}

\section*{Conflict of interest}

The author declares no conflict of interest.

\end{document}